\documentclass[11pt,final,a4paper]{amsart}
\usepackage{amssymb}
\usepackage{amsmath}
\usepackage{amsfonts}

\newcommand{\C}{\mathbb{C}}

\newcommand{\Q}{\mathbb{Q}}

\newcommand{\Z}{\mathbb{Z}}
\newcommand{\F}{\mathbb{F}}

\newtheorem{dummy}{Dummy}

\newtheorem{lemma}[dummy]{Lemma}
\newtheorem{theorem}[dummy]{Theorem}
\newtheorem{prop}[dummy]{Proposition}
\newtheorem{cor}[dummy]{Corollary}

\theoremstyle{definition}

\theoremstyle{remark}
\newtheorem{rem}[dummy]{Remark}

\begin{document}
\bibliographystyle{amsalpha}
\author{Sandro Mattarei}

\email{mattarei@science.unitn.it}

\urladdr{http://www-math.science.unitn.it/\~{ }mattarei/}

\address{Dipartimento di Matematica\\
  Universit\`a degli Studi di Trento\\
  via Sommarive 14\\
  I-38050 Povo (Trento)\\
  Italy}

\title{On a special congruence of Carlitz}

\begin{abstract}
We prove that if $q$ is a power of a prime $p$ and $p^{k}$ divides $a$, with $k\ge 0$, then
\[
1+(q-1)\sum_{0\le b(q-1)<a} \binom{a}{ b(q-1)}\equiv 0\pmod{p^{k+1}}.
\]
The special case of this congruence where $q=p$ was proved by Carlitz in 1953
by means of rather deep properties of the Bernoulli numbers.
A more direct approach produces our generalization and several related results.
\end{abstract}




\thanks{Partially supported by Ministero dell'Istruzione, dell'Universit\`a  e
  della  Ricerca, Italy,  through PRIN ``Graded Lie algebras  and pro-$p$-groups:
  representations, periodicity and derivations''.}

\maketitle

\thispagestyle{empty}

\section{Introduction}\label{sec:intro}

Sums of binomial coefficients of the form
\[
\sum_b \binom{a}{ b\,d+r}=
|\{X\subseteq\{1,\ldots,a\}:
|X|\equiv r\pmod{d}\}|
\]
occur in combinatorics and number theory.
Several classical results give information on the values of such sums modulo a prime or prime power.
One of the oldest results of this type is due to Hermite, who proved in 1876 that
(in modern notation) a prime $p$ divides
$\sum_{b>0} \binom{a}{ b(p-1)}$ if $a$ is a positive odd integer
(cf.~\cite[p.~271]{Dickson1}).
Hermite's result was then generalized in a number of directions,
the earliest due to Glaisher in 1899 (cf.~\cite[p.~272]{Dickson1}).
Glaisher showed that
\begin{equation}\label{eq:Glaisher}
\sum_{b\ge 0} \binom{a}{b(p-1)+r}\equiv
\binom{\bar a}{r}\pmod{p}
\end{equation}
for $p$ a prime, $a$ a positive integer, and $1\le r\le p-1$, where $\bar a$ denotes the smallest positive
integers congruent to $a$ modulo $p-1$.
This can also be formulated by saying that the value modulo $p$ of the left member of~\eqref{eq:Glaisher}
is a periodic function of $a>0$ with period $p-1$.
A proof of Glaisher's result based on Lucas' theorem for evaluating binomial coefficients
modulo a prime can be found in~\cite[Section~6]{Gra:organic}, but see the Introduction of~\cite{Sun:sums_binomial} for a simpler proof.
In Section~\ref{sec:multisection} we present an easy generalization of Glaisher's result which gives an efficient formula
for the value modulo $p$ of the sum
$\sum_b \binom{a}{b\,d+r}$, where $d$ is any integer prime to $p$.

In 1953 Carlitz~\cite{Carlitz:special} generalized Hermite's theorem to a prime power modulus by showing that
\begin{equation}\label{eq:Carlitz}
p+(p-1)\sum_{0<b(p-1)<a} \binom{a}{ b(p-1)}\equiv 0\pmod{p^{k+1}}
\end{equation}
if $p$ is an odd prime and $p^k$ divides the positive integer $a$.
(This is trivially true also for $p=2$, because the left member equals $2^a$ in this case.)
Unlike the proofs mentioned above of Glaisher's congruence~\eqref{eq:Glaisher},
Carlitz's proof of~\eqref{eq:Carlitz} is quite sophisticated.
It relies on certain congruences satisfied
by the Bernoulli numbers, namely that $B_m/m$ is a $p$-adic integer if $(p-1)\nmid m$, see~\cite[p.~238]{IR},
and that $(B_m+p^{-1}-1)/m$ is a $p$-adic integer if $(p-1)\mid m$, see~\cite[p.~247]{IR},
the latter result being due to Carlitz himself.
It seems that no other proof of Carlitz's congruence has ever appeared, except for the special case where $p-1$ divides $a$,
which follows from~\cite[Corollary~1.1]{SunTau}.

It appears most natural to prove Carlitz's congruence~\eqref{eq:Carlitz} by multisection of series.
In fact, this route allows us to prove the following generalization,
which does not seem amenable to Carlitz's original approach.

\begin{theorem}\label{thm:Carlitz}
If $q$ is a power of a prime $p$ and $p^{k}$ divides $a$, with $k\ge 0$, then
\[
1+(q-1)\sum_{0\le b(q-1)<a} \binom{a}{ b(q-1)}\equiv 0\pmod{p^{k+1}}.
\]
\end{theorem}

Although our approach does not allow an evaluation modulo $p^{k+1}$ (for $k>0$) of the more general sum
$\sum_b \binom{a}{ b(q-1)-r}$ where $r$ is any integer
(except for the case where $q-1$ divides $a$, considered in Corollary~\ref{cor:q-1} below),
it does produce the following remarkable symmetry.

\begin{theorem}\label{thm:symmetry}
Let $p$ be a prime, $q=p^f$, let $h,k$ be nonnegative integers with $h\ge k$ and $f\mid h+k$, and let $r,s$ be positive integers.
Then we have
\[
(-1)^{ps}\sum_{b} \binom{s\,p^k}{ b(q-1)-r}\equiv
(-1)^{pr}\sum_{b} \binom{r\,p^h}{ b(q-1)-s}
\pmod{p^{k+1}}.
\]
\end{theorem}

The case of Theorem~\ref{thm:symmetry} where $q-1$ divides $s$ has the following consequence,
which complements Theorem~\ref{thm:Carlitz} in the special case where $q-1$ divides $a$.

\begin{cor}\label{cor:q-1}
If $q$ is a power of a prime $p$, the number $(q-1)p^{k}$ divides $a$, with $k\ge 0$, and $q-1$ does not divide $r$, then
\[
(q-1)\sum_{b} \binom{a}{b(q-1)-r}\equiv -(-1)^{pr}\pmod{p^{k+1}}.
\]
\end{cor}

Carlitz went further in~\cite{Carlitz:special}
by evaluating the left member of his congruence~\eqref{eq:Carlitz}
modulo $p^{k+2}$, in terms of the Bernoulli numbers $B_{2s}$ and Wilson's quotient
$w_p=((p-1)!+1)/p$.
With notation slightly adapted to our present needs, Carlitz's result reads
\begin{multline}\label{eq:sharper}
p^{-k-1}\bigg\{1+(p-1)\sum_{0\le b(p-1)<sp^k} \binom{sp^k}{ b(p-1)}\bigg\}
\\
\equiv
s\bigg\{\frac{1}{2}
-\sum_{\substack{0<2j<sp^k\\ p-1\nmid 2j}} \binom{sp^k-1}{ 2j-1}\frac{B_{2j}}{2j}+\delta_{s}\frac{w_p}{p-1}\bigg\}
\pmod{p}
\end{multline}
for $p\ge 3$ (although stated for $p>3$ in~\cite{Carlitz:special}, see the beginning of our Section~\ref{sec:Carlitz_sharper}),
where $\delta_s=1$ if $p-1\mid s-1$ and $\delta_s=0$ otherwise.

Congruence~\eqref{eq:sharper} is useless for the purpose of a fast evaluation modulo $p$ of its left member,
because its right member is more complicated than the former, and contains more summands.
The first of our couple of contributions to~\eqref{eq:sharper}
is a proof that the value modulo $p$ of the left member of~\eqref{eq:sharper}
is actually independent of $k$, which is not apparent from the form of the right member.
In particular, the left member can be most conveniently evaluated modulo $p$ by replacing $k$ with $0$,
thus reducing the summation to about $s/(p-1)$ binomial coefficients.
Although the deeper connection with Bernoulli numbers shown by Carlitz's sharper congruence~\eqref{eq:sharper}
does not extend in an obvious way with a prime power $q$ replacing $p$,
we keep with the spirit of our previous results by allowing a prime power $q$ in place of $p$.

\begin{theorem}\label{thm:sharper}
Let $q=p^f$ be a power of a prime $p$, and let $s$ be a positive integer.
Then the value modulo $p$ of the expression
\[
p^{-k-1}\bigg\{1+(q-1)\sum_{0\le b(q-1)<sp^k} \binom{sp^k}{ b(q-1)}\bigg\},
\]
as a function of $k\ge 0$ (but $k\ge 2$ if $p=2$ and $s=1$),
depends only on the remainder of $k$ modulo $f$.
\end{theorem}

It is quite easy to see, in a way which we point out in Remark~\ref{rem:reduction},
that when $p$ is odd the value modulo $p$ of the expression considered in Theorem~\ref{thm:sharper}
is also a periodic function of $s>0$, with period dividing $(q-1)p$.
Thus, one only needs consider the range $0<s\le (q-1)p$.
Our final result displays a symmetry in the dependency on $s$ which allows one to further restrict this range in certain cases.
This is the only one among our results where we need to assume the prime $p$ to be odd.
We expand on the reasons for this after its proof in Section~\ref{sec:Carlitz_sharper}.

\begin{theorem}\label{thm:sharper_symmetry}
Let $q=p^f$ be a power of an odd prime $p$,
and let $k$ be a nonnegative integer.
Then the value modulo $p^{k+2}$ of the expression
\[
1+(q-1)\sum_{0\le b(q-1)<sp^k} \binom{sp^k}{ b(q-1)}
\]
is unaffected by replacing the positive integer $s$ with $p^t-s$,
where $t$ is any integer such that $p^t>s$ and $f\mid k+t$.
\end{theorem}

Computer calculations performed by means of the symbolic manipulation package {\tt MAPLE} have been extremely useful
for discovering and checking congruences, notably those of Theorems~\ref{thm:sharper} and~\ref{thm:sharper_symmetry}.

\section{Glaisher's congruence}
\label{sec:multisection}

A standard way of dealing with sums like that of Glaisher is based on
the identity
\begin{equation}\label{eq:multisection}
  \sum_b \binom{a}{ bd+r} x^{bd+r}
  =
  \frac{1}{d}
  \sum_{j=0}^{d-1}
  \omega^{-ir}
  (1+\omega^i x)^a
\end{equation}
in the polynomial ring $F[x]$,
where $F$ is any field containing a primitive $d$th root of unity $\omega$.
This is identity~(1.53) in~\cite{Gou}
(where $F$ is the complex field and $\omega=\exp(2\pi\mathrm{i}/d)$),
and follows by applying the more general formula for multisection of series~\cite[Chapter 1, Exercise 26]{Com}
to the generating function of the binomial coefficients,
$(1+x)^n=\sum_m\binom{n}{m}x^m$.
Here we follow the standard convention that unrestricted summation indices run over the integers;
however, the sum in~\eqref{eq:multisection} is a finite sum since $a$ is positive integer,
if $\binom{a}{ c}$ is defined to be $0$ for $c<0$, as usual.
The following result is the generalization of Glaisher's congruence announced in the Introduction.

\begin{prop}\label{prop:Glaisher}
Let $p$ be a prime, let $a$ and $d$ be positive integers with $p\nmid d$, and let $r$ be an integer.
If $f$ is the period of $p$ modulo $d$ and $\bar a$ is the smallest positive integer congruent to $a$ modulo $p^f-1$,
then we have
\[
\sum_b \binom{a}{ bd+r}
\equiv
\sum_b \binom{\bar a}{ bd+r}
\pmod{p}.
\]
In particular, for $d=q-1$ with $q$ a power of $p$ we have
\[
\sum_b \binom{a}{b(q-1)+r}
\equiv
\binom{\bar a}{ r}+
\binom{\bar a}{ r+q-1}
\pmod{p}.
\]

\end{prop}

\begin{proof}[First proof]
Let $\omega$ a primitive $d$th root of unity in the finite field of $q=p^f$ elements $\F_q$.
By evaluating the identity~\eqref{eq:multisection} for $x=1$ we obtain
\begin{equation*}
  \sum_b \binom{a}{ bd+r}
  =
\frac{1}{d}
  \sum_{\alpha\in \langle\omega\rangle}
  \alpha^{-r} (1+\alpha)^a.
\end{equation*}
The desired conclusion follows since the right member of the equality, as a function of the positive integer $a$,
depends only on the value of $a$ modulo $q-1$.
\end{proof}

We give another proof which does not use multisection of series.
\begin{proof}[Second proof]
The sum $\sum_b \binom{a}{bd+r}$ equals the coefficient of $x^r$ in the reduction of
$(1+x)^a$ modulo $x^d-1$.
If $d$ divides $q-1$, then $x^d-1$ divides $x^{q-1}-1$, and we have
\[
(1-x)^q=1-x^q\equiv 1-x\pmod{x^d-1}.
\]
Consequently, $(1-x)^{a+(q-1)}\equiv (1-x)^a\pmod{x^d-1}$ if $a>0$, and the conclusion follows.
\end{proof}

\begin{rem}\label{rem:d=q-1}
The general case of Proposition~\ref{prop:Glaisher} can also be deduced from its special case $d=q-1$,
by writing
$\sum_b \binom{a}{ bd+r}$
as
$\sum_{j=1}^{(q-1)/d}\sum_b \binom{a}{ b(q-1)+r+jd}$,
where $d$ divides $q-1$.
\end{rem}

\section{Carlitz's congruence}\label{sec:Carlitz}

Carlitz's congruence can be read as an equality in the ring $\Z/p^{k+1}\Z$.
We could then prove it by applying a version of identity~\eqref{eq:multisection} over this ring.
In fact, it is easily shown that identity~\eqref{eq:multisection} holds in the polynomial ring $R[x]$,
where $R$ is any commutative ring such that $d\cdot 1$ is invertible in $R$,
and $\omega$ is a unit of $R$
such that $\omega^d=1$ but $\omega^i-1$ is not a zero-divisor of $R$ for $0<i<d$.
Instead of this approach, we adopt here the equivalent but more standard way
of working in the (algebraic) integers and computing modulo $p^{k+1}$.
Nevertheless, a crucial ingredient of our proof of Carlitz's congruence would be
the following basic fact concerning the finite ring $\Z/p^{k+1}\Z$,
for $p$ odd~\cite[Chapter~4]{IR}:
its group of units is the direct product of two cyclic groups, one of order $p-1$ and one of order $p^k$.

In order to generalize Carlitz's congruence and prove Theorem~\ref{thm:Carlitz}
we need a corresponding result for the ring $R=\mathcal{O}/p^{k+1}\mathcal{O}$,
where $\omega$ is a primitive complex $(q-1)$-th root of unity and
$\mathcal{O}$ is the ring of integers in the cyclotomic field $\Q(\omega)$.
Note that $\mathcal{O}/p\mathcal{O}$ is the field of $q$ elements $\F_q$.

\begin{lemma}\label{lemma:units}
The group of units of $R=\mathcal{O}/p^{k+1}\mathcal{O}$, for $k\ge 0$,
is the direct product of a cyclic group of order $q-1$ and the group $1+pR$ of order $q^k$.
When $p$ is odd the latter is isomorphic with the additive group $pR$,
and hence has exponent $p^k$.
When $p=2$, the group $1+2R$ is the direct product of its subgroup $\{\pm 1\}$
and a subgroup isomorphic with a subgroup of index two of the additive group $2R$;
in particular, $1+2R$ has exponent $2^{k-1}$ if $q=2$, and $2^k$ if $q>2$.
\end{lemma}

One can prove Lemma~\ref{lemma:units} in an elementary way by induction on $k$ and similar calculations as those
performed in the standard proof, given in~\cite[Chapter~4, \S1]{IR}, of its special case where $q=p$
(namely, Equation~\eqref{eq:powers} in our proof of Theorem~\ref{thm:sharper} in the next section).
Such a proof can be found in~\cite{McDonald},
for example, where our Lemma~\ref{lemma:units} appears as~Theorem~XVI.9,
viewing $R$ as the Galois ring $GR(p^{k+1},f)$.
However, the following proof in the context of local fields seems more illuminating.

\begin{proof}
The finite quotient ring $R$ is unaffected if we replace $\Q(\omega)$ with its completion
$\Q_p(\omega)$ with respect to the prime divisor $p\mathcal{O}$.
In other words, we may work in the algebraic closure $\C_p$ of the field $\Q_p$ of $p$-adic numbers,
and let $\omega$ be a primitive $(q-1)$th root of unity in $\C_p$.
Then $K=\Q_p(\omega)$ is an unramified extension of $\Q_p$ of degree $f$, hence with residue field $\F_q$.
If
$O=\{x\in K:v_p(x)\ge 0\}$
is the valuation ring of $K$, and
$P=\{x\in K:v_p(x)\ge 1\}$
is the maximal ideal of $O$,
then $O/P^{k+1}\cong\mathcal{O}/p^{k+1}\mathcal{O}=R$.

According to~\cite[(III.4.4)]{Robert}, the group of units of $O$ splits into a direct product
$\mu_{q-1}\times(1+P)$,
where $\mu_{q-1}$ is the group of $(q-1)$th roots of unity in $\C_p$, which is generated by $\omega$.
Suppose first that $p$ is odd.
Since $1+P$ does not contain any nontrivial root of unity with $p$-power order, \cite[(V.4.2)]{Robert} shows that
the logarithm map
\[
1+\gamma\mapsto\sum_{j\ge 0}(-1)^{j-1}\gamma^j/j
\]
maps the multiplicative group $1+P$ isomorphically and isometrically onto the additive group $P$.
In particular, it maps $1+P^j$ onto $P^j$, for all positive integers $j$.
Consequently, the logarithm map induces an isomorphism
of the multiplicative group $(1+P)/(1+P^{k+1})$ onto the additive group $P/P^{k+1}\cong pR$.
However, $(1+P)/(1+P^{k+1})$ is the image of $1+P$ in the quotient ring $O/P^{k+1}$, and hence is isomorphic with $1+pR$, as claimed.

Suppose now that $p=2$, and assume that $k>0$ as we may.
Again according to~\cite[(V.4.2)]{Robert}, the logarithm map gives a group homomorphism of $1+P$ into $P$
with kernel $\mu_2=\{\pm 1\}$, which is the set of roots of unity of $2$-power order in $1+P$.
Its restriction to $1+P^2$ is an isometry onto $P^2$, and
hence maps $1+P^j$ bijectively onto $P^j$ for every $j\ge 2$.
Because the index $q$ of $1+P^2$ in $1+P$ equals the index of $P^2$ in $P$,
the logarithm maps $1+P$ onto a subgroup of index two of $P$.
Since, as before, $1+2R$ is isomorphic with $(1+P)/(1+P^{k+1})$,
its quotient $(1+2R)/\{\pm 1\}$ is isomorphic with a subgroup of index two of $2R$, call it $A$.
This leaves only two possibilities for the group structure of $1+2R$:
either it is isomorphic with $2R$, or it is the direct product of its subgroup $\{\pm 1\}$
and a subgroup isomorphic with $A$.
The former possibility would entail that $-1$, being an element of order two, should belong to
$(1+2R)^2\le 1+4R$, and is therefore to be excluded.
\end{proof}

The crucial part of Lemma~\ref{lemma:units} needed in the following proofs is the fact that
the exponent of $1+pR$ divides $p^k$.

\begin{proof}[Proof of Theorem~\ref{thm:Carlitz}]
Let $\omega$ be a primitive complex $(q-1)$-th root of unity and let
$\mathcal{O}$ be the ring of integers in the cyclotomic field $\Q(\omega)$.
According to identity~\eqref{eq:multisection} evaluated for $x=1$, in $\mathcal{O}$ we have
\[
1+(q-1)\sum_{0\le b(q-1)\le a} \binom{a}{ b(q-1)}=
\sum_{\alpha\in \mu_{q-1}\cup\{0\}}(1+\alpha)^a,
\]
where $\mu_{q-1}=\langle\omega\rangle$.
Since the elements $\alpha\in\mu_{q-1}\cup\{0\}$ are a set of representatives for the cosets
of the additive subgroup $p\mathcal{O}$ of $\mathcal{O}$, so are the elements $1+\alpha$.

Now view the above equality in the quotient ring
$R=\mathcal{O}/p^{k+1}\mathcal{O}$,
denoting by $\bar\mu_{q-1}$ the image of $\mu_{q-1}$ in $R$.
In particular, because $\beta(1+pR)=\beta+pR$ for $\beta\in R$,
the elements $1+\alpha$ for $\alpha\in\bar\mu_{q-1}\cup\{0\}\setminus\{-1\}$
are a set of representatives for the cosets of $1+pR$ in the group of units $U$ of $R$.
According to Lemma~\ref{lemma:units}, the group $U$ is the direct product of its subgroups
$\bar\mu_{q-1}$ and $1+pR$, and the latter has exponent $p^k$
(or $p^{k-1}$ when $q=2$, a trivial case here).
Consequently, if $\beta$ ranges over a set of representatives
for the cosets of $1+pR$ in $U$, then $\beta^{p^k}$
ranges over the elements of $\bar\mu_{q-1}$.
Taking into account also the case where $\alpha=-1$, it follows that
the elements $(1+\alpha)^{p^k}$ for $\alpha\in\bar\mu_{q-1}\cup\{0\}$
are distinct and coincide with the elements of $\bar\mu_{q-1}\cup\{0\}$.
Hence, in the ring $R$ we have
\begin{align*}
1+(q-1)\sum_{0\le b(q-1)\le a} \binom{a}{ b(q-1)}
&=
\sum_{\gamma\in\bar\mu_{q-1}\cup\{0\}}\gamma^{a/p^k}
\\&=
\begin{cases}
0 &\textrm{if $(q-1)\nmid a$},\\
q-1 &\textrm{if $(q-1)\mid a$}.
\end{cases}
\end{align*}
In both cases it follows that
\[
1+(q-1)\sum_{0\le b(q-1)< a} \binom{a}{ b(q-1)}=0
\]
in $R$, which is equivalent to the desired conclusion.
\end{proof}

In general, it does not seem possible to evaluate similarly modulo $p^{k+1}$ the more complicated right member
of the identity
\[
\sum_{b} \binom{a}{ b(q-1)-r}=
\frac{1}{q-1}\sum_{\alpha\in \mu_{q-1}}\alpha^{r}(1+\alpha)^a.
\]
However, one can somehow interchange the roles of the elements $\alpha$ and $1+\alpha$ in this formula,
as in the following proof.

\begin{proof}[Proof of Theorem~\ref{thm:symmetry}]
We adopt the same setting and notation as in the proof of Theorem~\ref{thm:Carlitz}.
We postpone to the end consideration of the case $p=2$ and assume first that $p$ is odd.

Let $\alpha\in\bar\mu_{q-1}\setminus\{-1\}$ and consider the element $\beta=(-1-\alpha)^{p^k}$ of $R$.
It is invertible and different from $-1$, because $\beta\equiv -1-\alpha^{p^k}\not\equiv 0,-1\pmod{\bar P}$,
where $\bar P$ denotes the image of $P$ in $R$, the unique maximal ideal of $R$.
Lemma~\ref{lemma:units} then implies that $\beta$
has multiplicative order dividing $q-1$, and hence belongs to $\bar\mu_{q-1}\setminus\{-1\}$.
Hence the correspondence $\alpha\mapsto\beta=(-1-\alpha)^{p^k}$ maps $\bar\mu_{q-1}\setminus\{-1\}$ into itself,
and so does the map $\beta\mapsto\alpha=(-1-\beta)^{p^h}$, because $h\ge k$.
We claim that these maps are inverse of each other.

In fact, if $\beta=(-1-\alpha)^{p^k}$ then
$(-1-\beta)^{p^h}\equiv(\alpha^{p^k})^{p^h}=\alpha\pmod{\bar P}$,
because $\alpha^q=\alpha$.
Since both $(-1-\beta)^{p^h}$ and $\alpha$ belong to $\bar\mu_{q-1}$, the congruence must be an equality.
Thus, the correspondence $\alpha\mapsto\beta=(-1-\alpha)^{p^k}$ is a permutation of $\bar\mu_{q-1}\setminus\{-1\}$.

Consequently, in $R$ we have
\begin{align*}
(-1)^{ps}\sum_{b} \binom{sp^k}{ b(q-1)-r}
&=
\frac{1}{q-1}\sum_{\alpha\in\bar\mu_{q-1}\setminus\{-1\}}\alpha^{r}(-1-\alpha)^{sp^k}
\\&=
\frac{1}{q-1}\sum_{\beta\in\bar\mu_{q-1}\setminus\{-1\}}(-1-\beta)^{rp^h}\beta^{s}
\\&=
(-1)^{pr}\sum_{b} \binom{rp^h}{ b(q-1)-s},
\end{align*}
which is equivalent to the desired conclusion.

The only difference in the case $p=2$ is that $\beta$ is not invertible when $\alpha=1$ and,
in fact, $\beta=(-1-1)^{2^k}=0$, because $2^k\ge k+1$.
Hence, in this case the map $\alpha\mapsto\beta=(-1-\alpha)^{p^k}$ does not send
$\bar\mu_{q-1}\setminus\{-1\}=\bar\mu_{q-1}$ into itself, but it does send
$\bar\mu_{q-1}\setminus\{1\}$ into itself.
Therefore, the final calculation remains valid by reading the summations over $\bar\mu_{q-1}\setminus\{1\}$
rather than over $\bar\mu_{q-1}\setminus\{-1\}$.
\end{proof}

\begin{rem}
When either $r=0$ or $s=0$, but not both, the congruence given in Theorem~\ref{thm:symmetry}
would be off by a summand $1/(q-1)$, as one can verify by going through the above proof in this anomalous situation.
In fact, this statement is equivalent to the special case of Theorem~\ref{thm:Carlitz} where
$p^k$ is a power of $q$.
\end{rem}

\begin{rem}
Because of the explicit formulas
$\sum_{k}\binom{n}{k}=2^n$ and
$\sum_{k}\binom{n}{2k}=\sum_{k}\binom{n}{2k+1}=2^{n-1}$,
which follow from Equation~\eqref{eq:multisection},
in the special cases where $q=2$ or $q=3$ the congruence stated in Theorem~\ref{thm:symmetry} reads
$2^{sp^k}\equiv 2^{rp^h}\pmod{2^{k+1}}$, and
\[
(-1)^s\cdot 2^{s\cdot 3^k-1}\equiv (-1)^r\cdot 2^{r\cdot 3^h-1}\pmod{3^{k+1}},
\]
which are easy to verify directly.
\end{rem}

\begin{rem}\label{rem:reduction}
Proposition~\ref{prop:Glaisher} (with our first proof, using Lemma~\ref{lemma:units})
extends at once to deal with congruences modulo a prime power $p^{k+1}$.
In view of Remark~\ref{rem:d=q-1}
this extension boils down to the following statement, for $p$ odd:
for integers $a$ and $r$ with $a$ a positive multiple of $p^k$ we have
\[
\sum_b \binom{a}{ b(q-1)+r}
\equiv
\sum_b \binom{\bar a}{ b(q-1)+r}
\pmod{p^{k+1}},
\]
where $\bar a$ is the smallest positive integer congruent to $a$ modulo $(q-1)p^k$.
(When $p=2$ the assertion holds only by taking $a,\bar a\ge k+1$,
because of the exceptional role of $\alpha=1$ in the proof of Theorem~\ref{thm:symmetry}.)
In the special case where $p^k$ divides $a$,
such an assertion can also be deduced (in an admittedly twisted way) from Theorem~\ref{thm:symmetry},
where the left member of the congruence is unaffected by adding to $s$ any multiple of $q-1$,
because the right member does.
\end{rem}

\begin{proof}[Proof of Corollary~\ref{cor:q-1}]
Write $a=sp^k$, thus $q-1\mid s$, and hence $(-1)^s=1$.
Choose an integer $h$ such that $h\ge k$ and $q-1\mid h+k$.
Using Theorem~\ref{thm:symmetry} and Theorem~\ref{thm:Carlitz} in turn we obtain
\begin{align*}
(q-1)\sum_{b} \binom{s\,p^k}{ b(q-1)-r}
&\equiv
(-1)^{pr}\cdot(q-1)\sum_{b} \binom{r\,p^h}{ b(q-1)}\equiv -(-1)^{pr},
\end{align*}
the congruences being modulo $p^{k+1}$.
\end{proof}

\begin{rem}
The above proof of Corollary~\ref{cor:q-1}
exploits the special case of Theorem~\ref{thm:symmetry} where $q-1$ divides $s$ but not $r$.
In contrast, the special case of Theorem~\ref{thm:symmetry} where both $r$ and $s$ are multiples of $q-1$,
which reads
\[
\sum_{b} \binom{s\,p^k}{ b(q-1)}\equiv
\sum_{b} \binom{r\,p^h}{ b(q-1)}
\pmod{p^{k+1}},
\]
yields no new information, since both sides of the congruences are already
known to be congruent to $q/(q-1)$ according to Theorem~\ref{thm:Carlitz}.
\end{rem}

\section{Carlitz's sharper congruence}\label{sec:Carlitz_sharper}

As we have pointed out in the Introduction, Carlitz stated congruence~\eqref{eq:sharper}
under the stronger hypothesis $p>3$.
In fact, his proof relies on a congruence for Bernoulli numbers which is valid only when $p>3$.
However, when $p=3$ congruence~\eqref{eq:sharper} remains valid by interpreting as zero the empty summation
in the right member.
In fact, the left member equals $(2^s-(-1)^s)/3^{k+1}$, which is easily seen to be congruent to
$s(1+\delta_{s})/2$ modulo $3$.
There appears to be no obvious interpretation of Equation~\eqref{eq:sharper} for $p=2$,
where its left member equals $2^{s2^k-k-1}$, and hence is congruent to $0$ modulo $2$
except when $s=1$ and $k=0$ or $1$.

Computer calculations show that the exceptional behaviour of~\eqref{eq:sharper} for small values of $s$ and $k$ when $p=2$
persists when generalizing to a power $q$ of $2$, as the statement of Theorem~\ref{thm:sharper} reflects.
The reason for this will be clear at the end of its proof.

\begin{proof}[Proof of Theorem~\ref{thm:sharper}]
Continue with the setting introduced in the proof of Theorem~\ref{thm:Carlitz}.
Thus, let $\omega$ be a primitive complex $(q-1)$-th root of unity, let
$\mathcal{O}$ be the ring of integers in the cyclotomic field $\Q(\omega)$,
and let $\mu_{q-1}=\langle\omega\rangle$.
We have
\[
1+(q-1)\sum_{0\le b(q-1)< sp^k} \binom{sp^k}{ b(q-1)}=
-\delta_s\cdot(q-1)+\sum_{\alpha\in \mu_{q-1}\cup\{0\}}(1+\alpha)^{sp^k},
\]
where $\delta_s=1$ if $q-1\mid s-1$ and $\delta_s=0$ otherwise.

We now assume that $p$ is odd and postpone a discussion of the case $p=2$ to the last paragraph of the proof.
For each $\alpha\in \mu_{q-1}\cup\{0\}$ we can write
\[
(1+\alpha)^{p^k}=\beta_\alpha(1+p^{k+1}\gamma_\alpha).
\]
with (uniquely determined) $\beta_\alpha\in\mu_{q-1}\cup\{0\}$ and $\gamma_\alpha\in\mathcal{O}$.
In fact, apart from the trivial case where $\alpha=-1$,
Lemma~\ref{lemma:units} implies that the image of $(1+\alpha)^{p^k}$
in the quotient ring $R=\mathcal{O}/p^{k+1}\mathcal{O}$
belongs to the image $\bar\mu_{q-1}$ of $\mu_{q-1}$.
Consequently, there exists a unique $\beta_\alpha\in\mu_{q-1}$ such that
$(1+\alpha)^{p^k}\beta_\alpha^{-1}\in 1+p^{k+1}\mathcal{O}$, as desired.

For every positive integer $n$ we have
\begin{equation}\label{eq:powers}
(1+t)^n\equiv 1+nt\pmod{pnt\mathcal{O}}
\end{equation}
provided $t\in p\mathcal{O}$ for $p$ odd, and $t\in 4\mathcal{O}$ for $p=2$.
This can be proved by extending standard calculations in the integers done in~\cite[Chapter~4, \S1]{IR},
but can also be deduced from its slightly more elegant $p$-adic version given in~\cite[(III.4.3)]{Robert}.
Since $p^{k+1}>2$ by hypothesis, it follows that
\[
(1+\alpha)^{sp^k}\equiv\beta_\alpha^s(1+sp^{k+1}\gamma_\alpha)
\pmod{p^{k+2}},
\]
and because $\beta_\alpha^q=\beta_\alpha$ we also have
\[
(1+\alpha)^{sp^kq}
\equiv
\beta_\alpha^s(1+sp^{k+1}q\gamma_\alpha)
\pmod{p^{k+2}q}.
\]

An argument seen in the Proof of Theorem~\ref{thm:Carlitz} shows that
$\beta_\alpha$ ranges over $\mu_{q-1}\cup\{0\}$ when $\alpha$ ranges over $\mu_{q-1}\cup\{0\}$, and hence
\[
-\delta_s\cdot(q-1)+\sum_{\alpha\in \mu_{q-1}\cup\{0\}}\beta_\alpha^s=0.
\]
Consequently, we have
\begin{multline*}
p^{-k-1}q^{-1}\bigg\{1+(p-1)\sum_{0\le b(p-1)<sp^kq} \binom{sp^kq}{ b(p-1)}\bigg\}
\equiv
\sum_{\alpha\in \mu_{q-1}\cup\{0\}}
\beta_\alpha^ss\gamma_\alpha
\\
\equiv
p^{-k-1}\bigg\{1+(p-1)\sum_{0\le b(p-1)<sp^k} \binom{sp^k}{ b(p-1)}\bigg\}
\pmod{p},
\end{multline*}
which implies the desired conclusion.

The peculiarity of the case $p=2$ is that
$(1+\alpha)^{2^k}$ cannot be expressed in the form $\beta_\alpha(1+2^{k+1}\gamma_\alpha)$ when $\alpha=1$.
However, this discrepancy has no consequences if we just set $\beta_1=0$, provided
$2^{s2^k}\equiv 0\pmod{2^{k+2}}$,
which is satisfied except when $s=1$ and $k=0$ or $1$.
\end{proof}

\begin{proof}[Proof of Theorem~\ref{thm:sharper_symmetry}]
With notation as in the Proof of Theorem~\ref{thm:sharper} we have
\[
1+(q-1)\sum_{0\le b(q-1)< sp^k} \binom{sp^k}{ b(q-1)}=
\sum_{\alpha\in \mu_{q-1}\setminus\{-1\}}\bigl((1+\alpha)^{sp^k}-\beta_\alpha^s\bigr),
\]
where $\beta_\alpha$ is the unique element of $\mu_{q-1}$ which is congruent to
$(1+\alpha)^{p^k}$ modulo $p$.
We already know from Theorem~\ref{thm:Carlitz} that the expression at the right member belongs to $p^{k+1}\mathcal{O}$.

For each $\alpha\in\mu_{q-1}\setminus\{-1\}$, let $\tilde\alpha$
be the unique element of $\mu_{q-1}\setminus\{-1\}$ which is congruent to $-\alpha/(1+\alpha)$ modulo $p$;
equivalently, let  $\tilde\alpha\in\mu_{q-1}\setminus\{-1\}$ be defined by the condition
$1+\tilde\alpha\equiv(1+\alpha)^{-1}\pmod{p}$.
The desired conclusion will follow from the congruence
\begin{equation}\label{eq:hard1}
(1+\tilde\alpha)^{sp^k}-\beta_{\tilde\alpha}^{s}
\equiv
(1+\alpha)^{(p^t-s)p^k}-\beta_\alpha^{p^t-s}
\pmod{p^{k+2}},
\end{equation}
which we prove in the following paragraphs.

Lemma~\ref{lemma:units} implies that
$\beta_\alpha\equiv(1+\alpha)^{p^kq}\pmod{p^{k+2}}$.
Note that we actually have
$\beta_\alpha\equiv(1+\alpha)^{p^kq^i}\pmod{p^{k+2}}$
for every $i>0$, because
$\gamma^{uq}=\gamma^u$ for all $\gamma\in\mathcal{O}/p^{k+2}\mathcal{O}$
if $p^{k+1}$ divides $u$.
We will use this observation without mention in the sequel
by multiplying or dividing certain exponents by appropriate powers of $q$ whenever convenient.
Since $\beta_{\tilde\alpha}=\beta_\alpha^{-1}$,
the claimed congruence~\eqref{eq:hard1} can be written in the equivalent form
\begin{multline}\label{eq:hard2}
(1+\tilde\alpha)^{sp^k}-(1+\alpha)^{-sp^kq}
\\
\equiv
(1+\alpha)^{(p^t-s)p^k}-(1+\alpha)^{(p^t-s)p^kq}
\pmod{p^{k+2}}.
\end{multline}

According to Lemma~\ref{lemma:units} we have
$\tilde\alpha\equiv
-\alpha/(1+\alpha)^{p^{k+t}}
\pmod{p^{k+2}}$.
Therefore, the left member of congruence~\eqref{eq:hard2} satisfies
\begin{align*}
(1+\tilde\alpha)^{sp^k}-&(1+\alpha)^{-sp^kq}
\equiv
\left(1-\frac{\alpha}{(1+\alpha)^{p^{k+t}}}\right)^{sp^k}-(1+\alpha)^{-sp^kq}
\\
&\equiv
(1+\alpha)^{-sp^kq}\left(\bigl((1+\alpha)^{p^{k+t}}-\alpha\bigr)^{sp^k}-1\right)
\\
&\equiv
(1+\alpha)^{-sp^kq}sp^k\bigl((1+\alpha)^{p^{k+t}}-\alpha-1\bigr)
\pmod{p^{k+2}}
\\
&=
(1+\alpha)^{1-sp^kq}sp^k\bigl((1+\alpha)^{p^{k+t}-1}-1\bigr)
\end{align*}
where we have used Equation~\eqref{eq:powers} in the next-to-last passage, because
$(1+\alpha)^{p^{k+t}}-\alpha\in 1+p\mathcal{O}$.
Similarly, the right member of congruence~\eqref{eq:hard2} satisfies
\begin{align*}
(1+\alpha)^{(p^t-s)p^k}-&(1+\alpha)^{(p^t-s)p^kq}
\equiv
(1+\alpha)^{p^{k+t}-sp^k}-(1+\alpha)^{p^{k+t}-sp^kp^{k+t}}
\\
&\equiv
(1+\alpha)^{p^{k+t}-sp^kp^{k+t}}
\bigl((1+\alpha)^{(p^{k+t}-1)sp^k}-1\bigr)
\\
&\equiv
(1+\alpha)^{p^{k+t}-sp^kq}
sp^k
\bigl((1+\alpha)^{p^{k+t}-1}-1\bigr)
\pmod{p^{k+2}}
\end{align*}
where we have used Equation~\eqref{eq:powers} in the last passage, this time because
$(1+\alpha)^{p^{k+t}-1}\in 1+p\mathcal{O}$.

In order to complete a proof of congruence~\eqref{eq:hard2}, and hence of Theorem~\ref{thm:sharper_symmetry},
it suffices to observe that
\[
(1+\alpha)^{p^{k+t}-1}-1\equiv
(1+\alpha)^{p^{k+t}-1}\bigl((1+\alpha)^{p^{k+t}-1}-1\bigr)
\pmod{p^2}.
\]
In fact, this congruence is equivalent to
\[
\bigl((1+\alpha)^{p^{k+t}-1}-1\bigr)^2\equiv 0\pmod{p^2},
\]
which holds because
$(1+\alpha)^{p^{k+t}-1}-1\equiv 0\pmod{p}$.
\end{proof}

The above proof breaks down for $p=2$, mainly because of the two appeals to Equation~\eqref{eq:powers},
which requires $t\in 4\mathcal{O}$ rather than $t\in 2\mathcal{O}$ when $p=2$.
In fact, computer calculations show that the statement of~Theorem~\ref{thm:sharper_symmetry}
fails for $p=2$ and $q>4$, even subject to (reasonable) restrictions on $k$.
We have considered the trivial case where $q=2$ earlier in this section.
When $q=4$, one can show by means of Equation~\eqref{eq:multisection} that
the expression considered in Theorem~\ref{thm:sharper_symmetry} equals $2^{s2^k}$ for $k>0$
and is, therefore, a multiple of $2^{k+2}$
except when $(s,k)=(1,1)$.

\bibliography{References}

\end{document}